\documentclass{amsart}

\usepackage{amsmath,amssymb}

\newtheorem{theorem}{Theorem}
\newcommand{\bt}{\begin{theorem}}
\newcommand{\et}{\end{theorem}} 

\newtheorem{corollary}{Corollary}
\newcommand{\bc}{\begin{corollary}}
\newcommand{\ec}{\end{corollary}} 

\newtheorem{lemma}{Lemma}
\newcommand{\bl}{\begin{lemma}}
\newcommand{\el}{\end{lemma}}

\newtheorem{problem}{Problem}
\newcommand{\bprob}{\begin{problem}}
\newcommand{\eprob}{\end{problem}}

\newcommand{\beq}{\begin{equation}}
\newcommand{\eeq}{\end{equation}}

\newcommand{\A}{\ensuremath{ \mathbf A }}
\newcommand{\C}{\ensuremath{ \mathbf C }}
\newcommand{\F}{\ensuremath{ \mathbf F }}
\newcommand{\N}{\ensuremath{ \mathbf N }}

\newcommand{\R}{\ensuremath{\mathbf R}}
\newcommand{\benum}{\begin{enumerate}}
\newcommand{\eenum}{\end{enumerate}}

\newcommand{\mcm}{\mathcal{M}}

\DeclareMathOperator{\qqand}{\qquad\text{and}\qquad}
\DeclareMathOperator{\card}{\text{card}}

\title[$B_h$-sets and perturbations]{$B_h$-sets and perturbations in normed vector spaces}

\author{Melvyn B. Nathanson}

\address{Department of Mathematics, Lehman College (CUNY), Bronx, NY, USA}  
\email{melvyn.nathanson@lehman.cuny.edu}

\date{\today}


\begin{document}

\begin{abstract}
The subset $A = \{a_i:i \in I\}$ of a normed vector space is a $B_h$-set if every element 
of the sumset $hA$ has a unique representation as a sum of $h$ elements of $A$.  
Let $\varepsilon = \{\varepsilon_i:i \in I\}$ be a set of positive real numbers.  
An $\varepsilon$-perturbation of $A$ is a set $A' = \{a'_i:i\in I\}$ such that $|a_i'-a_i|<\varepsilon_i$ 
for all $i \in I$. 
Let $\Delta_{hA} = \inf\{|x'-x| : x,x' \in hA \text{ and } x\neq x'\}$.  
It is proved that if $A$ is  finite or countably infinite set with $\Delta_{hA}>0$, 
then there is a $B_h$-set $A'$ that is an $\varepsilon$-perturbation of $A$. 
\end{abstract}

\maketitle 


Let $\N = \{1,2,3,\ldots\}$ be the set of positive integers 
and  let $I = \N$ or let $I = \{1,2,\ldots, k\}$ for some $k \in \N$. 
For $h \in \N$, let $\{m_i:i \in I\}$ be a set of nonnegative integers 
such that 
\[
\sum_{i\in I} m_i = h. 
\]
Let $\mcm(h,I)$ be the set of all such sets.  

Let $X$ be an additive abelian group or semigroup and let $A = \{a_i:i \in I\}$ be a  subset of $X$ 
with $a_i \neq a_j$ for all $i \neq j$. 
The \emph{$h$-fold sumset }of $A$ is the set 
\[
hA = \left\{ \sum_{i\in I} m_i a_i:  \{m_i:i \in I\} \in \mcm(h,I) \right\}.
\]
Thus, $hA$ is the set of all sums of $h$ not necssarily distinct elements of $A$.  

For all $x \in X$, the \emph{representation function} $r_{A,h}(x)$ counts the number of 
representations of $x$ as a sum of $h$ elements of $A$, that is,
\[
r_{A,h}(x) = \card\left\{   \{m_i:i \in I\}  \in \mcm(h,I):  \sum_{i \in I} m_ia_i = x\right\}.
\]
The set $A$ is a \emph{$B_h$-set} if $r_{A,h}(x) = 1$ for all $x \in hA$, that is, 
if every element in $X$ has at most one representation as a sum of $h$ elements of $A$. 
 A $B_2$-set is also called a \emph{Sidon set}. 

Almost all additive and combinatorial number theory work on $B_h$-sets has been 
on $B_h$-sets in the integers or in other discrete groups.  
Cilleruelo  and Ruzsa~\cite{cill-ruzs04} and Nathanson~\cite{nath21,nath26} 
have studied continuous analogues for $B_h$-sets of real, complex, and $p$-adic numbers.  
In this note we consider $B_h$-sets in normed  vector spaces over a field \F\ of  characteristic 
of 0.

Let $X \neq \{0\}$ be a normed vector space over the field $\F$ 
and let $\varepsilon = \{\varepsilon_i:i \in I\}$ be a set of positive real numbers.  
Let $A = \{a_i:i \in I\}$ and $A' = \{a'_i:i \in I\}$ be subsets of $X$. 
The set $A'$ is an \emph{$\varepsilon$-perturbation} of $A$ if 
\[
|a'_i - a_i | < \varepsilon_i 
\]
for all $i \in I$.

\bl       \label{perturbe:lemma:Bh}
Let $X \neq \{0\}$ be a normed  vector space over $\F$ 
and let $\varepsilon = \{\varepsilon_i:i \in I\}$ be a set of positive real numbers. 
For every integer $h \geq 2$, the space $X$ contains uncountably many  
$B_h$-sets $W = \{w_i :i\in I\}$ such that 
\[
| w_i| < \varepsilon_i 
\]
for all $i \in I$. 
\el

\begin{proof} 
Choose an integer $g \geq h+1$.  The uniqueness of the $g$-adic 
representation of a positive rational number implies that the set 
\[
\left\{ \frac{1}{g^i} : i \in I \right\}
\] 
is a $B_h$-set.  
There are uncountably many  strictly increasing sequences $(j_i)_{i=0}^{\infty}$ 
of positive integers 
such that 
\[
\frac{1}{g^{j_i}} < \varepsilon_i
\]
for all $i \in I$.  Then 
\[
\left\{ \frac{1}{g^{j_i}} : i \in I \right\}
\] 
is a $B_h$-set. 
For every vector $w \in X$ with $|w| = 1$, the set 
\[
\left\{ w_i =  \frac{1}{g^{j_i}} w: i \in I \right\}
\] 
is a $B_h$-set with $|w_i| < \varepsilon_i$ for all $i \in I$.  
This completes the proof. 
\end{proof}

For every subset $A$ of a  normed vector space $X$ with $|A| > 1$, 
we define 
\[
\Delta_{A} = \inf \{|a-a'|:a,a' \in A  \text{ and } a \neq a'\}. 
\] 
If $A$ is finite, then $\Delta_A> 0$.
Let $\varepsilon = \{\varepsilon_i:i \in I\}$ be a set of positive real numbers. 
A $\varepsilon$-perturbation of a subset $A = \{a_i:i \in I\}$ of $X$ is a set 
$A' = \{a'_i:i \in I\}$ in $X$ such that $|a'_i - a_i| < \varepsilon_i$ for all $i \in I$.

\bt                         \label{perturbe:theorem:real}   
Let  $h \geq 2$ and let $A = \{a_i:i \in I\}$ be a subset of the normed  vector space $X$ 
such that $\Delta_{hA} > 0$.  
For every set $\varepsilon = \{\varepsilon_i:i \in I\}$ of positive real numbers, 
there are uncountably many  $B_h$-sets $A' = \{a'_i:i \in I\}$ in $X$ that are 
$\varepsilon$-perturbations of $A$. 
\et

\begin{proof} 
For all $i \in I$, choose $\varepsilon'_i$ such that 
\[
0 < \varepsilon_i' <\min\left( \varepsilon_i, \frac{\Delta_{hA}}{2h} \right) 
\] 
and let $\varepsilon' = \{\varepsilon'_i:i \in I\}$. 
By Lemma~\ref{perturbe:lemma:Bh}, there are uncountably many 
$B_h$-sets $\{w_i : i \in I\}$ in $X$ such that  $|w_i |< \varepsilon'_i$ 
for all $i \in I$.   For each such set, let 
\[
A' = \{ a'_i = a_i + w_i: i \in I \}.
\]
Then 
\[
|a'_i-a_i| = |w_i| < \varepsilon'_i < \varepsilon_i 
\]
and so $A'$ is an $\varepsilon$-perturbation of $A$.  
We shall prove that $A'$ is a $B_h$-set in $X$. 

Let $\{m_i:i \in I\}$ and $\{m'_i:i \in I\}$ be distinct sets in $\mcm(h,I)$ 
and let $ \sum_{i \in I} m_ia'_i$ and $ \sum_{i \in I} m'_ia'_i$ be the 
associated elements of the sumset $hA'$. 
Then 
\[
 \sum_{i \in I} m'_ia'_i = \sum_{i \in I} m_ia'_i
\]
if and only if 
\[
 \sum_{i \in I} m_i(a_i +w_i)= \sum_{i \in I} m'_i (a_i + w_i)
\]
if and only if 
\[
 \sum_{i \in I} m_i a_i  -  \sum_{i \in I} m'_i a_i 
=   \sum_{i \in I} m_i w_i  -  \sum_{i \in I} m'_i w_i.
\]
We have $ \sum_{i \in I} m_i a_i  \in hA$ and $ \sum_{i \in I} m'_i a_i  \in hA$.  
If $ \sum_{i \in I} m_ia_i  \neq \sum_{i \in I} m'_i a_i $, then 
\begin{align*}
\Delta_{hA} & \leq \left|  \sum_{i \in I} m_ia_i - \sum_{i \in I} m'_i a_i \right| 
= \left|  \sum_{i \in I} m_i w_i  -  \sum_{i \in I} m'_i w_i\right| \\ 
&  \leq\sum_{i \in I} m_i  \left|   w_i \right|  + \sum_{i \in I} m'_i  \left|  w_i\right| 
 \leq  \sum_{i \in I} ( m_i + m'_i ) \varepsilon'_i  \\ 
& <2h \left( \frac{\Delta_{hA}}{2h} \right) = \Delta_{hA}
\end{align*}
which is absurd.  

If $ \sum_{i \in I} m_ia_i = \sum_{i \in I} m'_i a_i $, then 
\[
 \sum_{i \in I} m_i w_i  = \sum_{i \in I} m'_i w_i. 
\] 
Because $\{w_i : i \in I\}$ is a $B_h$-set, 
it follows that $\{m_i:i \in I\} \neq \{m'_i:i \in I\}$,  
which is also absurd. Therefore, $A'$ is a $B_h$-set. 
This completes the proof. 
\end{proof}

\end{document}